%% file: gaps.tex
\magnification 1050
 \def\oui{oui}
 \ifx\textures\oui
 \input CrayolaColors 
 \long\def\rge#1{\Red#1\Black}

 \else 
 \input itenn.tex

\beginpackages
\usepackage{color}
\endpackages 
\long\def\rge#1{{\color{red}#1}}

 \fi
\newcount\numchapi\numchapi=0

\input option_keys.tex
\optionparag=2
\def\paradouze{oui}
\anglais
\dimstand\hautspages{G. Tenenbaum}{Normal largest gap between prime factors}\dateheure

\titrecentre{A note on the normal largest  gap between prime factors}
\bigskip\medskip
\centerline{GŽrald Tenenbaum}
\bigskip\bigskip
{\eightpoint\leftskip1cm\rightskip1cm
\noi{\bf Abstract.} Let $\{p_j(n)\}_{j=1}^{\omega(n)}$ denote the increasing sequence of distinct prime factors of an integer~$n$. We provide details for the proof of a statement of Erd\H{o}s implying that, for any function $\xi(n)$ tending to infinity with $n$, we have
 $$f(n):=\max_{1\leqslant j<\omega(n)}\log \Big({\log p_{j+1}(n)\over \log p_j(n)}\Big)=\log_3n+O(\xi(n))$$  
 for almost all integers $n$.
\PAR
\medskip\noi
 \par
\smallskip 
\noi \bf 2010 Mathematics Subject Classification: \rm  11N25, 11N35.\par }
\par 
\bigskip\bigskip
In private correspondence \citer{So18}, E. Sofos asked for normal upper and lower bounds for the arithmetic function
$$f(n):=\max_{1\leqslant j<\omega(n)}\log \Big({\log p_{j+1}(n)\over \log p_j(n)}\Big),$$
where $\{p_j(n)\}_{j=1}^{\omega(n)}$ denotes the increasing sequence of distinct prime factors of an integer $n$.
The answer to this question is actually provided by a statement of Erd\H os in \citer{Er59}, where it is asserted without proof that, for any $c>0$,
the density of the set of those integers $n$ such that $\e^{f(n)}>c\log_2n$ is $1-\e^{-1/c}$. Here and in the sequel $\log_k$ denotes the $k$-fold iterated logarithm.
Indeed, given any function $\xi(n)$ tending to infinity and letting $c$ tend sufficiently slowly to zero, we infer that, for almost all $n$, we have $f(n)>\log_3n-\xi(n)$, while by letting $c$ tend slowly to infinity we obtain that $f(n)\leqslant \log_3n+\xi(n)$ also holds for a set of asymptotic density~1. Interesting related results appear in Erd\H os' articles \citer{Er66} and \citer{Er69}.
\par 
In this short note, we provide the not so obvious details of the proof, in the spirit of the paper \citer{Er59}.\par 
\propt{encfn}{ (Erd\H os \citer{Er59})}{Let $c>0$. Then the inequality $f(n)>\log_3n-\log (1/c)$ holds on a set of integers $n$ of asymptotic density $1-\e^{-1/c}$. In particular, given any function $\xi(n)\to\infty$,  we have
$$\log_3n-\xi(n)\leqslant  f(n)\leqslant\log_3n+\xi(n)$$
for almost all integers $n$.}
\nid Let $x$ be a large parameter, write $Z:=c\log_2x$, and let $\chi_p(n)$ denote the indicator function of the set of those integers $n\leqslant x$ 
that are divisible by the prime $p$ but by no prime $q$ such that $p< q\leqslant  p^Z$. For squarefree $m$, put $\chi_m(n):=\prod_{p\mid m}\chi_p(n)$. Taking into account the slow growth of $\log_2x$, it is sufficient to show that 
$$N(x):=|\{n\leqslant x:p\mid n,\,p\leqslant x^{1/Z}\Rightarrow \chi_p(n)=0\}|\sim \e^{-1/c}x\qquad (x\to\infty).\eqdef{defN}$$
Indeed, $n\leqslant x$ is counted by $N(x)$ if, and only if, $p_{j+1}(n)\leqslant p_j(n)^Z$ whenever $1\leqslant j<\omega(n)$.  Let $P^+(n)$ (resp. $P^-(n))$ denote the largest (resp. the smallest) prime factor of an integer $n$, with the convention that $P^+(1)=1$ (resp. $P^-(1)=\infty$), and let $\mu$ designate the Mšbius function. By the inclusion--exclusion principle (see, \eg, \citer{Te15}, \page39) or the Mšbius inversion formula, we thus have, for $x\geqslant 1$,
$$N(x)=\sum_{\di{m\leqslant x}{m\in\M}}\mu(m)\sum_{n\leqslant x}\chi_m(n),\eqdef{fii}$$
where $\M$ denotes the set of all squarefree integers $m\leqslant x$ such that $P^+(m)\leqslant x^{1/Z}$ and $p_{j+1}(m)>p_j(m)^Z$ whenever $j\in[1,\omega(m)[$.
Moreover, as is familiar in Brun's sieve method, we obtain a lower or an upper  bound for the left-hand side by restricting the outer sum to integers $m$ with  number of prime factors at most equal to an odd or an even bound.
\par \goodbreak
Let $k\geqslant 1$ be fixed. For each $m\in\M$, define $P_m:=\prod_{p|m}\prod_{p< q\leqslant p^Z}q$. Then
$$N_k(x):=\sum_{\di{m\in\M}{\omega(m)=k}}\sum_{n\leqslant x}\chi_m(n)=\sum_{\di{m\in\M}{\omega(m)=k}}\sum_{\di{\nu\leqslant x/m}{(\nu,P_m)=1}}1.\eqdef{Nk}$$
  From a standard sieve result (see, \eg, the lemma in \citer{Te93}), we infer that the inner sum~is 
 $$\sim{x\varphi(P_m)\over mP_m}={x\over mZ^k}\Big\{1+O\Big({1\over \log P^-(m)}\Big)+o(1)\Big\}\eqdef{inn}$$
 provided $P^+(m)\leqslant x^{o(1/Z)}$. The contribution of the remaining integers $m$ will be treated as an error term. We may plainly assume that $x^{\varepsilon_x}<P^+(m)^Z\leqslant x$, with $\varepsilon_x$ tending to 0 arbitrarily  slowly. Then
the inner sum is classically (see, \eg, \citer{Te15}, Exercise 85) $$\ll{x\varphi(P_m)\over mP_m}\ll{x\over mZ^k}
$$
and the corresponding contribution to $N_k(x)$ is
$${x\over Z^k}\sum_{x^{\varepsilon_x/Z}<p\leqslant x^{1/Z}}{1\over p}\sum_{\di{P^+(h)\leqslant x}{\omega(h)=k-1}}{1\over h}\ll{x\log (1/\varepsilon_x)\over \log_2x}=o(x).$$
Inserting \eqref{inn} back into \eqref{Nk} and summing over $m$, we obtain, for $k\geqslant 0$,
$$N_k(x)={x\over Z^k}\sum_{\di{p_k\leqslant x^{1/Z}}{p_{j}\leqslant p_{j+1}^{1/Z}\,(1\leqslant j<k)}}{1+O(1/\log p_1)+o(1)\over p_1\cdots p_k}\sim{x(\log_2x)^k\over k!Z^k}\sim {x\over c^kk!}\cdot$$ 
  Thus, for arbirary $\ell\geqslant 1$, we have
$$x\sum_{0\leqslant k\leqslant 2\ell+1}{(-1)^k\over c^kk!}+o(x)\leqslant N(x)\leqslant x\sum_{0\leqslant k\leqslant 2\ell}{(-1)^k\over c^kk!}+o(x),$$
and the required result follows by selecting $\ell$ arbitrarily large. 
 \qed\bigskip\goodbreak

\medskip
\noi{\bf Acknowledgements.} The author wishes to express warm thanks to RŽgis de la Bretche for valuables  discussions on this matter.

\bigskip
\centerline{\twelvebf References}\bigskip
\eightpoint{
\bibtem{Er59} P. Erd\H os,
Some remarks on prime factors of integers.
{\it Canad. J. Math. \bf11} (1959), 161--167. \par 
\bibtem{Er66} P. Erd\H os, 
On some properties of prime factors of integers,
{\it Nagoya Math. J. \bf27} (1966), 617--623.\par 
\bibtem{Er69} P. Erd\H os, On the distribution of prime divisors, \it
Aequationes Math. \bf2 \rm(1969), 177Ð183. \par 
\bibtem{So18} E. Sofos, private e-mail message, August 30, 2018.\par 
\bibtem{Te93} G. Tenenbaum, Cribler les entiers sans grand facteur premier, in: R. C. Vaughan (ed.), Theory
and applications of numbers without large prime factors, {\it Phil. Trans.  R.  Soc. Lond.} A {\bf
345} (1993), 377--384.\par 
\bibtem{Te15} G. Tenenbaum, {\it Introduction to analytic and probabilistic number theory}, 3rd ed., Graduate Studies in Mathematics 163, Amer. Math. Soc. 2015.\par
\par }
\vskip 5mm
{\sevenrm\baselineskip9pt
G\'erald Tenenbaum\par
Institut \'Elie Cartan\par 
Universit\'e de Lorraine\par
 BP 70239\par
54506 Vand\oe uvre-ls-Nancy Cedex\par
 France
\smallskip
 \seventt gerald.tenenbaum@univ-lorraine.fr\par}
\bye

%% file: itenn.tex

  
  %
  %

  %
  %

  \catcode`@=12 

 \def\defrefnote#1{\definexref{#1}{{\the\footnotenumber}}{refnotes}}

  %
  %


 \input eplain.tex
\makeatletter
\def\numberedfootnote{%
ÊÊ\global\advance\footnotenumber by 1
ÊÊ\@eplainfootnote{{\number\footnotenumber}}%
}%
\def\makecolumns#1/#2 {\par \begingroup
ÊÊ \@columndepth = #1
ÊÊ \advance\@columndepth by -1
ÊÊ \divide \@columndepth by #2
ÊÊ \advance\@columndepth by 1
ÊÊ \@linestogoincolumn = \@columndepth
ÊÊ \@linestogo = #1
ÊÊ \currentcolumn = 1
ÊÊ \def\@endcolumnactions{%
ÊÊÊÊÊÊ\ifnum \@linestogo<2
ÊÊÊÊÊÊÊÊ \the\crtok \egroup \endgroup \par 
ÊÊÊÊÊÊ\else
ÊÊÊÊÊÊÊÊ \global\advance\@linestogo by -1
ÊÊÊÊÊÊÊÊ \ifnum\@linestogoincolumn<2
ÊÊÊÊÊÊÊÊÊÊÊÊ\global\advance\currentcolumn by 1
ÊÊÊÊÊÊÊÊÊÊÊÊ\global\@linestogoincolumn = \@columndepth
ÊÊÊÊÊÊÊÊÊÊÊÊ\the\crtok
ÊÊÊÊÊÊÊÊ \else
ÊÊÊÊÊÊÊÊÊÊÊÊ&\global\advance\@linestogoincolumn by -1
ÊÊÊÊÊÊÊÊ \fi
ÊÊÊÊÊÊ\fi
ÊÊ }%
ÊÊ \makeactive\^^M
ÊÊ \letreturn \@endcolumnactions
ÊÊ \@columnwidth = \hsize
ÊÊÊÊ \advance\@columnwidth by -\parindent
ÊÊÊÊ \divide\@columnwidth by #2
ÊÊ \penalty\abovecolumnspenalty
ÊÊ \noindent 
ÊÊ \valign\bgroup
ÊÊÊÊ &\hbox to \@columnwidth{\strut \hsize = \@columnwidth ##\hfil}\cr
}%
\makeatother

\lefteqnumbers
   \def\testd{oui}
   \def\choixlat{\ifx\numadroite\testd\righteqnumbers
            \else  \lefteqnumbers\fi}
    \choixlat

\catcode`@=\letter
\def\@eplainfootnote#1{\let\@sf\empty 
  \ifhmode\edef\@sf{\spacefactor\the\spacefactor}\/\fi
  \global\advance\hlfootlabelnumber by 1
  \hlstart@impl{foot}{\hlfootlabel}%
  \hldest@impl{footback}{\hlfootbacklabel}%
  \hbox{$^{(#1)}$}%
  \hlend@impl{foot}%
  \@sf\vfootnote{#1.}%
}%
\catcode`@=\other

  \interfootnoteskip=0pt
  
  \everyfootnote={\eightpoint\leftskip=5truemm\rightskip5truemm}
  
  \hsize150truemm\vsize 240truemm\hoffset=5truemm
  \def\dimstand{\hsize 150truemm\vsize 240truemm\hoffset=5truemm\voffset=0truemm}
  
  \pretolerance=500\tolerance=1000\brokenpenalty=5000
  \parindent3mm
  
  \countdef\temps=170
  \temps=\time
  \countdef\nminutes=171{\nminutes = \time}
  \countdef\nheure=172
  \def\heure{\begingroup                     
     \temps = \time \divide\temps by 60
     \nheure = \temps                        
     \nminutes = \time
     \multiply\temps by 60
     \advance\nminutes by -\temps            
     \ifnum\nminutes<10 \toks1 = {0}%
     \else\toks1 = {}%
     \fi
     \number\nheure h\the\toks1 \number\nminutes  
  \endgroup}%

  \newcount\chstart
  \chstart=\pageno
 \headline={\ifnum\pageno=\chstart {\hfill} \else {\hss \tenrm --\ \folio\ --\hss}\fi}
  \footline={\hfill}
  \normalbaselines
  \frenchspacing
    \def\dater{\vglue-10mm\rightline{(\the\day/\the\month/\the\year)}}
  \def\dateheure{\vglue-10mm\rightline{(\the\day/\the\month/\the\year,\ \heure)}}

  \newif\ifpagetitre \pagetitretrue
\newtoks\hautpagetitre \hautpagetitre={\hfill}
\newtoks\baspagetitre \baspagetitre={\hfill}
\newtoks\auteurcourant \auteurcourant={\hfill}
\newtoks\titrecourant \titrecourant={\hfill}
\newtoks\hautpagegauche
\newtoks\hautpagedroite
\newtoks\hautpagemilieu
\hautpagemilieu={\tenrm\hfil -- \folio\ -- \hfil}
\hautpagegauche={\ifx\midfolio\oui\the\hautpagemilieu\else\tenrm\folio\hfill\the\auteurcourant\hfill\fi}
\hautpagedroite={\ifx\midfolio\oui\the\hautpagemilieu\else\hfill\the\titrecourant\hfill\tenrm\folio\fi}
\newtoks\baspagegauche \baspagegauche={\hfil}
\newtoks\baspagedroite \baspagedroite={\hfil}
\headline={\ifpagetitre\the\hautpagetitre
\else\ifodd\pageno\the\hautpagedroite\else\the\hautpagegauche\fi\fi }
\footline={\ifpagetitre\the\baspagetitre
\else\ifodd\pageno\the\baspagedroite
\else\the\baspagegauche\fi\fi \global\pagetitrefalse}

\def\pageblanche{\vfill\eject\pagetitretrue
\null\vfill\eject
\pagetitretrue
}
\def\chgtpage{\ifodd\pageno \else
\pageblanche \fi \pagetitretrue\titreun=0\footnotenumber=0}

\def\chgtpageincrtitreun{\ifodd\pageno \else
\pageblanche \fi \pagetitretrue\footnotenumber=0}

\def\majnombres{\ifodd\pageno \else
\pageblanche \fi \pagetitretrue\hautpoly\titreun=0\footnotenumber=0}

\def\hautspages#1#2{\auteurcourant={\ninepcap#1}\titrecourant={\nineit#2}}

\ifnum\chstart=\pageno \pagetitretrue\fi
  

  \def\PAR{\par}
  
  \def\leftnote#1{\vadjust{\setbox1=\vtop{\hsize 20mm\parindent=0pt\eightpoint
  \baselineskip=9pt\rightskip=4mm plus 4mm\vskip-4.7mm#1}\hbox{\kern-2cm\smash{\box1}}}}

  
  \def\raggedcenter{\leftskip=20pt plus 10em  
       \rightskip=\leftskip 
        \parfillskip=0pt 
         \spaceskip=.3333em \xspaceskip=.5em 
          \pretolerance=9999 \tolerance=9999
           \hyphenpenalty=9999 \exhyphenpenalty=9999 }
           
  \def\titrecentre#1{{\parindent0mm\raggedcenter
       \spaceskip=.6em plus .2em minus .2em\xspaceskip=.6em plus .2em minus .2em
        \tit#1\par}}
        


  \def\oui{oui}
  
   \def\fontetitreun{\ifx\paradouze\oui\douzepts\gpdouze\twelvebf\textfont1=\twelveib\else
\quatorzepts\gpquatorze\fourteenbf\fi}

\def\fontetitreunl{\douzepts\textfont1=\twelveib\scriptfont1=\tenib\fourteenti}
 
 \def\fontetitredeux{\textfont1=\eleveni\ifx\paradouze\oui\onzepts\scriptfont1=\ninei\elevenit\else
                        \douzepts\twelveit\fi}
 
   \def\fontetitredeuxb{\ifx\paradouze\oui\onzepts\eleventi\gponze\textfont1=\elevenib\scriptfont1=\nineib
                         \else\douzepts\twelveti\scriptfont1=\twelveib\scriptfont1=\tenib\gpdouze\fi}
                         
\def\fontetitredeuxl{\onzepts\textfont1=\elevenbf\scriptfont1=\ninebf\twelvebf}
  
\def\fontetitretrois{\textfont0=\elevenrm\scriptfont0=\eightrm\textfont1=\eleveni
                      \scriptfont1=\eighti\scriptscriptfont1=\sixi\elevenit}
                      
\def\fontetitrequatre{\textfont0=\elevenrm\scriptfont0=\eightrm\textfont1=\eleveni
                      \scriptfont1=\eighti\scriptscriptfont1=\sixi\elevenrm}
  
  \newcount\titreun\titreun=0
  \newcount\titredeux\titredeux=0
  \newcount\titretrois\titretrois=0
  \newcount\titrequatre\titrequatre=0
  \newcount\enonce\enonce=0
  
  \def\incr#1{\global\advance#1 by 1 {\the #1}}
  \def\avance#1{\global\advance#1 by 1}
  \def\init#1{\global#1=0}
  
  \long\def\Indentation#1#2{\setbox10=\hbox{\fontetitreun#1}
  	                    \ifdim\wd10 < 4mm
                         \setbox10=\hbox to 4mm{\box10\hfill}
                       \else\ifdim\wd10 < 6mm
                         \setbox10=\hbox to 6mm{\box10\hfill}
  	                    \else\ifdim\wd10 < 8mm
                         \setbox10=\hbox to 8mm{\box10\hfill}
                       \else\ifdim\wd10 < 12mm
                         \setbox10=\hbox to 12mm{\box10\hfill}
                       \fi\fi\fi\fi
                       \dimen10=\hsize
                       \advance \dimen10 by -\wd10
                       \noindent \box10 %
                       \ignorespaces
                       \hbox{\vtop{\hsize=\dimen10\raggedright\noindent\fontetitreun#2}}}

  \long\def\paraun#1{\removelastskip\par\medskip\goodbreak\vskip0pt plus.01\vsize\penalty-100
                \vskip0pt plus-.01\vsize
  	              \init{\titredeux}\ifnum\optionparag=1{\init\eqnumber\init\enonce}\else{}\fi
                  \goodbreak{\fontetitreun
  	                \Indentation{\incr{\titreun}.\ }{\fontetitreun #1\par}}\nobreak\medskip}

 %
 %
 \long\def\paraunc#1{\removelastskip\par\bigskip\goodbreak\vskip0pt plus.01\vsize\penalty-100
                \vskip0pt plus-.01\vsize
  	              \init{\titredeux}
                 \ifnum\optionparag=1{\init{\eqnumber}\init\enonce}\else{}\fi
                  \goodbreak
  	                {\parindent0mm\raggedcenter\fontetitreun\incr{\titreun}.\ 
                     \fontetitreun #1\par}\nobreak\medskip}
                     
\newtoks\titreunl
\titreunl={\ifnum\titreun=1{I}\fi%
\ifnum\titreun=2{II}\fi%
\ifnum\titreun=3{III}\fi%
\ifnum\titreun=4{IV}\fi%
\ifnum\titreun=5{V}\fi%
\ifnum\titreun=6{VI}\fi%
\ifnum\titreun=7{VII}\fi%
\ifnum\titreun=8{VIII}\fi%
\ifnum\titreun=9{IX}\fi%
\ifnum\titreun=10{X}\fi%
\ifnum\titreun=11{XI}\fi%
\ifnum\titreun=12{XII}\fi%
\ifnum\titreun=13{XIII}\fi%
}
\long\def\paraunl#1{\removelastskip\par\bigskip\bigskip\goodbreak\vskip0pt plus.01\vsize\penalty-100
                \vskip0pt plus-.01\vsize
  	              \init{\titredeux}\ifnum\optionparag=1{\init\eqnumber\init\enonce}\else{}\fi
                  \goodbreak{\fontetitreunl
  	                \Indentation{\global\advance\titreun by 1{\the\titreunl}.\ }{\fontetitreunl #1\par}}\nobreak\smallskip}

  
  \long\def\paradeux#1{\init{\titretrois}\vskip0pt plus.01\vsize\penalty-10
                \vskip0pt plus-.01\vsize\ifx \elie\oui\medskip\ifnum\titredeux>0\medskip\fi\fi
                 \Indentation{\fontetitredeux\the\titreun${\cdot}$\incr{\titredeux}.}
                              {\fontetitredeux\textfont1=\eleveni#1}\nobreak\par }
  
  \long\def\paradeuxb#1{\init{\titretrois}\vskip0pt plus.001\vsize\penalty-10
                \vskip0pt plus-.01\vsize{\ifx \elie\oui\medskip\ifnum\titredeux>0\medskip\fi\fi
                  \Indentation
  {\fontetitredeuxb\the\titreun${\cdot}$\incr{\titredeux}.}{ \fontetitredeuxb#1}}\nobreak
\smallskip}

\newtoks\titredeuxl
\titredeuxl={\ifnum\titredeux=1{A}\fi%
\ifnum\titredeux=2{B}\fi%
\ifnum\titredeux=3{C}\fi%
\ifnum\titredeux=4{D}\fi%
\ifnum\titredeux=5{E}\fi%
\ifnum\titredeux=6{F}\fi%
\ifnum\titredeux=7{G}\fi%
\ifnum\titredeux=8{H}\fi%
\ifnum\titredeux=9{I}\fi%
\ifnum\titredeux=10{J}\fi%
\ifnum\titredeux=11{K}\fi%
\ifnum\titredeux=12{L}\fi%
\ifnum\titredeux=13{M}\fi%
}
 \long\def\paradeuxl#1{\init{\titretrois}\vskip0pt plus.001\vsize\penalty-10
                \vskip0pt plus-.01
                \vsize \bigskip%
                  \Indentation
     {\fontetitredeuxl\global\advance\titredeux by 1
  \quad \the\titreunl${\cdot}$\the\titredeuxl.}{ \fontetitredeuxl#1}
  \removelastskip\nobreak\smallskip}
  

  \long\def\paratrois#1{\init{\titrequatre}\ifdim\lastskip<\smallskipamount
                \removelastskip\smallskip\fi
                 \vskip0pt plus.01\vsize\penalty-10
                  \vskip0pt
plus-.01\vsize{\ifx \elie\oui\ifnum\titretrois>0\medskip\fi\fi
\Indentation{\fontetitretrois\the\titreun${\cdot}$\the\titredeux${\cdot}$\incr{\titretrois}.\ }
  {\hskip0mm\baselineskip=14pt\fontetitretrois#1}\nobreak\smallskip}}
  
  
  \long\def\paratroisl#1{\init{\titrequatre}\ifdim\lastskip<\smallskipamount
                \removelastskip\fi
                 \vskip0pt plus.01\vsize\penalty-10
                  \vskip0pt
plus-.01\vsize\ifx \elie\oui\bigskip
\fi
\Indentation{\fontetitretrois\quad \quad \the\titreunl{${\cdot}$}\the\titredeuxl${\cdot}$\incr{\titretrois}.\ }
  {\hskip0mm\fontetitretrois#1}\nobreak\smallskip}


  \long\def\paraquatre#1{\ifdim\lastskip<\smallskipamount
                \removelastskip\smallskip\fi
                 \vskip0pt plus.01\vsize\penalty-10
                  \vskip0pt
                  plus-.01\vsize\par
 
\Indentation{\fontetitrequatre \the\titreun{${\cdot}$}\the\titredeux${\cdot}$\the\titretrois${\cdot}$\incr{\titrequatre}.\ }
{\hskip0mm\fontetitrequatre#1}\nobreak\smallskip}


\newtoks\titrequatrel
\titrequatrel={\ifnum\titrequatre=1{a}\fi%
\ifnum\titrequatre=2{b}\fi%
\ifnum\titrequatre=3{c}\fi%
\ifnum\titrequatre=4{d}\fi%
\ifnum\titrequatre=5{e}\fi%
\ifnum\titrequatre=6{f}\fi%
\ifnum\titrequatre=7{g}\fi%
\ifnum\titrequatre=8{h}\fi%
\ifnum\titrequatre=9{i}\fi%
\ifnum\titrequatre=10{j}\fi%
\ifnum\titrequatre=11{k}\fi%
\ifnum\titrequatre=12{l}\fi%
\ifnum\titrequatre=13{m}\fi%
}
\long\def\paraquatrel#1{\ifdim\lastskip<\smallskipamount
                \removelastskip\smallskip\fi
                 \vskip0pt plus.01\vsize\penalty-10
                  \vskip0pt
                  plus-.01\vsize{\bigskip
\Indentation{\global\advance\titrequatre by 1
\fontetitrequatre\quad \quad \quad \the\titreunl${\cdot}$\the\titredeuxl${\cdot}$\the\titretrois${\cdot}$\the\titrequatrel.\ }
{\hskip0mm\fontetitrequatre#1}\nobreak\smallskip}}

\ifx\optionkeys\oui
\def\drefun#1{\definexref{¤#1}{{\the\titreun}}{}} 
\def\drefdeux#1{\definexref{¤#1}{{\the\titreun}.{\the\titredeux}}{}}
\def\dreftrois#1{\definexref{¤#1}{{\the\titreun}.{\the\titredeux}.{\the\titretrois}}{}}
\else
\def\drefun#1{\definexref{prg#1}{{\the\titreun}}{}} 
\def\drefdeux#1{\definexref{prg#1}{{\the\titreun}.{\the\titredeux}}{}}
\def\dreftrois#1{\definexref{prg#1}{{\the\titreun}.{\the\titredeux}.{\the\titretrois}}{}}
\fi

%


  \long\def\propdeux#1#2#3#4{%
       \avance{\enonce}
       \leavevmode\edef\temp{#2}%
         \ifx\temp\empty 
          \else
           \definexref{#2}{#1~{\the\titreun.\the\enonce}}{enonces}
            \definexref{s#2}{{\the\titreun.\the\enonce}}{enonces}
             \fi
\smallskip
      \noindent{\bf#1\ {\bf\the\titreun.\the\enonce{#3}.}\enspace}{\sl#4\par}%
      \ifdim\lastskip<\medskipamount \removelastskip\penalty55\par \fi
   }

  \long\def\propun#1#2#3#4{%
      \avance{\enonce}
       \leavevmode\edef\temp{#2}%
        \ifx\temp\empty 
          \else
           \definexref{#2}{#1~{\the\enonce}}{enonces}
            \definexref{{s#2}}{{\the\enonce}}{enonces}
             \fi
   \par 
     \noindent{\bf#1\ {\bf\the\enonce{#3}.}\enspace}{\sl#4\par}%
     \ifdim\lastskip<\medskipamount \removelastskip\penalty55\medskip\fi
  }
  
  \long\def\prop#1#2#3#4{\ifnum\optionparag=1
                          \propdeux{#1}{#2}{\textfont1=\elevenib#3}{#4} \else\propun{#1}{#2}{\textfont1=\elevenib#3}{#4}\fi}

  \long\def\propt#1#2#3{\ifx\tpf\oui \prop{Th\'eo\-r\`eme}{#1}{#2}{#3}\par
                       \else\prop{Theorem}{#1}{#2}{#3}\par\fi}
  
  \long\def\propl#1#2#3{\ifx\tpf\oui\prop{Lem\-me}{#1}{#2}{#3}\par
                         \else\prop{Lemma}{#1}{#2}{#3}\par\fi}
  
  \long\def\propc#1#2#3{\ifx\tpf\oui\prop{Corol\-laire}{#1}{#2}{#3}\par
                         \else\prop{Corollary}{#1}{#2}{#3}\par\fi}

  \long\def\propd#1#2#3{\ifx\tpf\oui\prop{D\'efi\-nition}{#1}{#2}{#3}\par
                       \else\prop{Definition}{#1}{#2}{#3}\par\fi} 
  
  \long\def\proptd#1#2#3{\ifx\tpf\oui\prop{Th\'eor\`eme et d\'efi\-nition}{#1}{#2}{#3}\par
                       \else\prop{Theorem and definition}{#1}{#2}{#3}\par\fi}


  
  \newcount\optionparag\optionparag=1
  
  \long\def\section#1#2{\ifnum\optionparag=1 \paraun{#2} 
                        \else\goodbreak{\fontetitreun
  	                \Indentation{#1.\ }{#2}}\nobreak\smallskip\fi}

  \def\eqconstruct#1{\ifnum\optionparag=1{\the\titreun\hbox{$\cdot$}#1}\else{#1}\fi}

  
  
  \def\numref{oui}  
  
  \newcount\mesref\mesref=0 
  \def\defbib#1{\ifx\numref\oui\global\advance\mesref by 1 \definexref{#1}{{\the
                 \mesref}}{}\else\definexref{#1}{#1}{}\fi}
  \def\bibtem#1{\defbib{#1}\item{\citer{#1}}}
  \def\citer#1{[\ref{#1}]}

  
  \font\seventeenmsa=msam10 at 17pt    
  \font\fourteenmsa=msam10 at 14pt
  \font\twelvemsa=msam10 at 12pt
  \font\tenmsa=msam10                 
  \font\ninemsa=msam10 at 9pt 
  \font\eightmsa=msam10 at 8pt 
  \font\sevenmsa=msam7 
  \font\sixmsa=msam10 at 6pt
  \font\fivemsa=msam5
  \newfam\msafam\textfont\msafam=\tenmsa\scriptfont\msafam=\sevenmsa\scriptscriptfont\msafam=\fivemsa
  
  \font\seventeenbb=msbm10 at 17pt     
  \font\fourteenbb=msbm10 at 14pt
  \font\twelvebb=msbm10 at 12pt
  \font\tenbb=msbm10                   
  \font\ninebb=msbm10 at 9pt 
  \font\eightbb=msbm10 at 8pt 
  \font\sevenbb=msbm7 
  \font\sixbb=msbm10 at 6pt
  \font\fivebb=msbm5 
  \newfam\bbfam\textfont\bbfam=\tenbb\scriptfont\bbfam=\sevenbb\scriptscriptfont\bbfam=\fivebb
  \def\bb{\fam\bbfam\tenbb}%

  \font\seventeenscaln=eusm10 at 17pt   
  \font\twelvescaln=eusm10 at 12pt
  \font\tenscaln=eusm10                
  \font\ninescaln=eusm10 scaled 900
  \font\eightscaln=eusm10 scaled 800
  \font\sevenscaln=eusm10 scaled 700
  \font\sixscaln=eusm10 scaled 600
   
  \newfam\scalnfam\textfont\scalnfam=\tenscaln\scriptfont\scalnfam=\sevenscaln\scriptscriptfont\scalnfam=\sixscaln
  \def\scaln{\fam\scalnfam\tenscaln}%
  \def\scal{\scaln}
  
  \font\tenscalb=eusb10                

  \font\sevenscalb=eusb10 scaled 700

  \newfam\scalbfam\textfont\scalbfam=\tenscalb\scriptfont\scalbfam=\sevenscalb
  %
  
  %
  %
  \font\fourteenrm=cmr12 scaled 1200
  \font\elevenrm=cmr10 at 11pt
  \font\twelverm=cmr12
  \font\ninerm=cmr9
  \font\eightrm=cmr8      
  \font\sevenrm=cmr7
  \font\sixrm=cmr6

  \font\seventeenpcap=cmcsc10 at 17pt
  \font\tenpcap=cmcsc10                        
  \font\ninepcap=cmcsc9
  \font\eightpcap=cmcsc8
  \font\sevenpcap=cmcsc10 scaled 700
  
  \newfam\pcapfam\textfont\pcapfam=\tenpcap\scriptfont\pcapfam=\sevenpcap
  \def\pcap{\fam\pcapfam\tenpcap}
  
  \font\seventeenrm=cmbx12 scaled 1400

  \font\fourteenbf=cmbx10 scaled 1400
  
  \font\twelvebf=cmbx12
  \font\elevenbf=cmbx10 at 11pt
  \font\ninebf=cmbx9  
  \font\eightbf=cmbx8
  \font\sixbf=cmbx6
  
  \font\tengot=eufm10                           
   
  \font\eightgot=eufm10 at 8truept 
  \font\sevengot=eufm7 
  \font\sixgot=eufm10 at 6 truept 
   
  \newfam\gotfam
  \textfont\gotfam=\tengot\scriptfont\gotfam=\sevengot\scriptscriptfont\gotfam=\sixgot
  %

  
  \def\tit{%
  \textfont0=\seventeenrm\scriptfont0=\tenrm\def\rm{\fam0\seventeenrm}%
  \textfont1=\seventeenib\scriptfont1=\twelveib%
  \textfont2=\seventeensy\scriptfont2=\twelvesy\scriptscriptfont2=\ninesy
  \textfont3=\seventeenex
  \textfont\itfam=\seventeenti
  \def\it{\fam\itfam\seventeenti}%
  \textfont\bbfam=\seventeenbb \scriptfont\bbfam=\twelvebb
  \def\bb{\fam\bbfam\seventeenbb}%
  \textfont\msafam=\seventeenmsa\scriptfont\msafam=\twelvemsa
  \textfont\scalnfam=\seventeenscaln
  \def\pcap{\fam\pcapfam\seventeenpcap}
  \normalbaselineskip=25pt\normalbaselines\rm}

  \font\seventeenti=cmbxti10 scaled 1680
  
  \font\fourteenti=cmbxti10 at 14pt
  
  \font\twelveti=cmbxti10 scaled 1200
  \font\eleventi=cmbxti10 at 11pt

  %
  %
  \font\twelveit=cmti12	
  \font\elevenit=cmti10 scaled 1100
  \font\nineit=cmti9
  \font\eightit=cmti8
  \font\sevenit=cmti7

  %
  %
 
 \font\seventeenib=cmmib10 scaled 1680
  \font\fourteenib=cmmib10 scaled 1400
  \font\twelveib=cmmib10 scaled 1200
  \font\elevenib=cmmib10 scaled 1100
  \font\tenib=cmmib10
\font\eightib=cmmib10 scaled 800
  \font\nineib=cmmib10 scaled 900
\font\sevenib=cmmib10 scaled 700
\font\sixib=cmmib10 scaled 600
\font\fiveib=cmmib10 scaled 500

\ifx\ITAN\oui
\else
\innernewfam\cmmibfam
\textfont\cmmibfam=\tenib
\scriptfont\cmmibfam=\sevenib
\scriptscriptfont\cmmibfam=\fiveib
\def\ib{\fam\cmmibfam\tenib}
\fi

  %
  %
  
  \font\eleveni=cmmi10 scaled 1100
  \font\ninei=cmmi9
  \font\eighti=cmmi8 
  \font\seveni=cmmi7 			                
  \font\sixi=cmmi6
  
  \font\ninesl=cmsl9                    
  \font\eightsl=cmsl8 
  \font\sevensl=cmsl10 at 7pt

  \font\ninett=cmtt9                    
  \font\eighttt=cmtt8
  \font\seventt=cmtt10 scaled 700

  \font\seventeensy=cmsy10 scaled 1680    
  \font\fourteensy=cmsy10 scaled 1400
  \font\twelvesy=cmsy10 scaled 1176
  
  \font\ninesy=cmsy9                      
  \font\eightsy=cmsy8
  \font\sixsy=cmsy6
  \font\seventeenex=cmex10 at 17pt
  \font\fourteenex=cmex10 at 14pt
  \font\twelveex=cmex10 at 12pt
  \font\nineex=cmex10 at 9pt
  \font\eightex=cmex10 at 8pt
  \font\sevenex=cmex10 at 7pt
  \font\sixex=cmex10 at 6pt
  \font\fiveex=cmex10 at 5pt
  
   
  \font\fourteengp=cmmi10 at 14pt
  
  \font\twelvegp=cmmib10 at 12pt
  \font\elevengp=cmmib10 at 11pt
  \font\tengp=cmmib10                          
  \font\ninegp=cmmib10 at 9pt 
  \font\eightgp=cmmib8 
   
  \font\sixgp=cmmib6


  \def\gponze{\textfont0=\elevenbf\scriptfont0=\eightbf\scriptscriptfont0=\sixbf
           \textfont1=\elevengp\scriptfont1=\eightgp\scriptscriptfont1=\sixgp}
  \def\gpdouze{\textfont0=\twelvebf\scriptfont0=\tenbf\scriptscriptfont0=\ninebf
           \textfont1=\twelvegp\scriptfont1=\tengp\scriptscriptfont1=\ninegp}        
  
 \def\gpquatorze{\textfont0=\fourteenbf\scriptfont0=\twelvebf\scriptscriptfont0=\elevenbf
           \textfont1=\fourteengp\scriptfont1=\twelvegp\scriptscriptfont1=\elevengp}

  
  \expandafter\chardef\csname pre amssym.def at\endcsname=\the\catcode`\@
  \catcode`\@=11
  \def\undefine#1{\let#1\undefined}
  \def\newsymbol#1#2#3#4#5{\let\next@\relax
   \ifnum#2=\@ne\let\next@\msafam@\else
   \ifnum#2=\tw@\let\next@\bbfam@\fi\fi
   \mathchardef#1="#3\next@#4#5}
  \def\mathhexbox@#1#2#3{\relax
   \ifmmode\mathpalette{}{\m@th\mathchar"#1#2#3}%
   \else\leavevmode\hbox{$\m@th\mathchar"#1#2#3$}\fi}
  \def\hexnumber@#1{\ifcase#1 0\or 1\or 2\or 3\or 4\or 5\or 6\or 7\or 8\or
   9\or A\or B\or C\or D\or E\or F\fi}
  
  \def\setboxz@h{\setbox\z@\hbox}
  \def\wdz@{\wd\z@}
  \def\boxz@{\box\z@}
  
  \edef\msafam@{\hexnumber@\msafam}
  \mathchardef\dabar@"0\msafam@39
  
  \edef\bbfam@{\hexnumber@\bbfam}
  \def\widehat#1{\setboxz@h{$\m@th#1$}%
   \ifdim\wdz@>\tw@ em\mathaccent"0\bbfam@5B{#1}%
   \else\mathaccent"0362{#1}\fi}
  \def\widetilde#1{\setboxz@h{$\m@th#1$}%
   \ifdim\wdz@>\tw@ em\mathaccent"0\bbfam@5D{#1}%
   \else\mathaccent"0365{#1}\fi}
  \newsymbol\leqq 1335          
  \newsymbol\leqslant 1336
  \newsymbol\lessgtr 1337       
  \newsymbol\backprime 1038     
  \newsymbol\risingdotseq 133A  
  \newsymbol\fallingdotseq 133B 
  \newsymbol\succcurlyeq 133C   
  \newsymbol\geqq 133D          
  \newsymbol\geqslant 133E
  \newsymbol\nmid 232D
  \newsymbol\nexists 2040
  \newsymbol\smallsetminus 2272
  \newsymbol\varnothing 203F
  
  \catcode`\@=\active

  \catcode`\@=11
  \newcount\typofr\typofr=1
  
  \catcode`\;=\active
  \def;{\ifnum\typofr=1\relax\ifhmode\ifdim\lastskip>\z@\unskip\fi
     \kern.2em\fi\string;\else\string;\fi}
  
  \catcode`\:=\active
  \def:{\ifnum\typofr=1\relax\ifhmode\ifdim\lastskip>\z@\unskip\fi
  \penalty\@M\ \fi\string:\else\string:\fi}
  
  \catcode`\!=\active
  \def!{\ifnum\typofr=1\relax\ifhmode\ifdim\lastskip>\z@\unskip\fi
     \kern.2em\fi\string!\else\string!\fi}
  
  \catcode`\?=\active
  \def?{\ifnum\typofr=1\relax\ifhmode\ifdim\lastskip>\z@\unskip\fi
     \kern.2em\fi\string?\else\string?\fi}

  \def\francais{\typofr=1\def\tpf{oui}}
  \def\anglais{\typofr=2\def\tpf{non}\def\english{oui}}
  \def\oui{oui}
  \francais
  
  \catcode`\@=12
  

\ifx\textures\oui
\def\raye #1|{\leavevmode\setbox1=\hbox{#1}%
\raise .5pt\hbox to \wd1{\xleaders\hbox{\rge{$ \sct / $}%
\kern 1pt}\hfill\kern -1pt }\kern -\wd1 \unhbox1\relax }

\def\barre #1|{\leavevmode\setbox1=\hbox{#1}%
\rlap{\Red\vrule height 2.4pt depth -1.2pt width \wd1}\Black \unhbox1\relax}
\else
\def\raye #1|{\leavevmode\setbox1=\hbox{#1}%
\raise .5pt\hbox to \wd1{\xleaders\hbox{\rge{$ \sct / $}%
\kern 1pt}\hfill\kern -1pt }\kern -\wd1 \unhbox1\relax }

\def\barre #1|{\leavevmode\setbox1=\hbox{#1}%
\rlap{\color{red}\vrule height 2.4pt depth -1.2pt width \wd1}\color{black} \unhbox1\relax}

\fi
  

  
  \def\og{\leavevmode\raise.24ex\hbox{$\scriptscriptstyle\langle\!\langle\>$}}    
  \def\fg{\leavevmode\raise.24ex\hbox{$\scriptscriptstyle\>\rangle\!\rangle$}}    

  \def\M{{\scal M}}
  
  \def\O{{\scal O}}
  \def\P{{\scaln P}}

  \def\frac#1#2{{#1\over #2}}
  \def\di#1#2{\sct#1\atop{\sct#2}}

  \def\page#1{\rm p.\thinspace#1}

  \def\qedbox{$\rlap{$\sqcap$}\sqcup$}           
  \def\qed{\nobreak\hfill\penalty250 \hbox{}\nobreak\hfill\qedbox\par }

  \def\¤{\S\thinspace}

  \def\¥{$\bullet$ }
  
  
  \def\e{{\rm e}}

  \def\epsilon{\varepsilon}

  \def\phi{\varphi}
  \def\theta{\vartheta}
  \def\rho{\varrho}

  \def\sct{\scriptstyle}
  \def\pf{\noi{\it Proof. }}
  \def\nid{\ifnum\typofr=1\par\noindent{\it D\'emonstration. }\else\pf\fi}
  \def\noi{\noindent}
  \def\rem{\ifnum\typofr=1\noi{\it Remarque.}\ \else\noi{\it Remark.}\ \fi}
  \def\rems{\ifnum\typofr=1\noi{\it Remarques.}\ \else\noi{\it Remarks.}\ \fi}

  \def\1{{\bf 1}}
  \def\|{\Vert}

  \def\leq{\leqslant}
  \def\geq{\geqslant}

  \def\eg{{e.g.}}
  

  \def\log{\mathop{\rm log}\nolimits}




  \def\pmb#1{\setbox0=\hbox{#1}%
  \kern-.025em\copy0\kern-\wd0\kern.05em\copy0\kern-\wd0\kern-.025em\raise .0433em\box0 }

  
  \skewchar\eighti='177 \skewchar\sixi='177
  \skewchar\eightsy='60 \skewchar\sixsy='60
  
  \def\eightpoint{%
  \textfont0=\eightrm\scriptfont0=\sixrm\scriptscriptfont0=\fiverm
  \def\rm{\fam0\eightrm}%
  \textfont1=\eighti\scriptfont1=\sixi
  \scriptscriptfont1=\fivei\def\oldstyle{\fam1\seveni}%
  \textfont2=\eightsy\scriptfont2=\sixsy\scriptscriptfont2=\fivesy
  \textfont3=\eightex\scriptfont3=\sixex
  \textfont\itfam=\eightit
  \def\it{\fam\itfam\eightit}%
  \textfont\slfam=\eightsl
  \def\sl{\fam\slfam\eightsl}%
  \textfont\bbfam=\eightbb \scriptfont\bbfam=\sixbb\scriptscriptfont\bbfam=\fivebb
  \def\bb{\fam\bbfam\eightbb}%
  \textfont\msafam=\eightmsa\scriptfont\msafam=\sixmsa
  \textfont\scalnfam=\eightscaln
  \def\scaln{\fam\scalnfam\eightscaln}
  \textfont\ttfam=\eighttt
  \def\tt{\fam\ttfam\eighttt}%
\textfont\gotfam=\eightgot
  \textfont\bffam=\eightbf\scriptfont\bffam=\sixbf\scriptscriptfont\bffam=\fivebf
  \def\bf{\fam\bffam\eightbf}%
  \ifx\ITAN\oui\else\textfont\cmmibfam=\eightib
       \scriptfont\cmmibfam=\sixib
        \scriptscriptfont\cmmibfam=\fiveib
         \def\ib{\fam\cmmibfam\eightib}
   \fi
  \textfont\pcapfam=\eightpcap
  \def\pcap{\fam\pcapfam\eightpcap}
  \abovedisplayskip=2pt plus2pt minus 2pt
  \belowdisplayskip=2pt plus1pt minus 2pt
  \abovedisplayshortskip= 1pt plus 2pt minus 1pt
  \belowdisplayshortskip= 1pt plus 2pt minus 1pt
  \smallskipamount=2pt plus 1pt minus 2pt
  \medskipamount=3pt plus 2pt minus 2pt
  \bigskipamount=7pt plus 3pt minus 3pt
  \setbox\strutbox=\hbox{\vrule height 5pt depth 2pt width 0pt}%
  \normalbaselineskip=9pt\normalbaselines\rm}

  \def\({\left(}
  \def\){\right)}
  
  \def\footnoterule{\kern -2pt\hrule width 7truecm\kern 2.4pt}
  
  \def\xnotedef#1{\definexref{#1}{\noexpand\number\footnotenumber}{Note}}%

  
  
  \def\ninepoint{%
  \textfont0=\ninerm\scriptfont0=\sixrm\scriptscriptfont0=\fiverm
  \def\rm{\fam0\ninerm}%
  \textfont1=\ninei\scriptfont1=\sixi
  \scriptscriptfont1=\fivei\def\oldstyle{\fam1\ninei}%
  \textfont2=\ninesy\scriptfont2=\sixsy\scriptscriptfont2=\fivesy
  \textfont3=\nineex\scriptfont3=\sixex
  \textfont\itfam=\nineit
  \def\it{\fam\itfam\nineit}%
  \textfont\slfam=\ninesl
  \def\sl{\fam\slfam\ninesl}%
  \textfont\bbfam=\ninebb\scriptfont\bbfam=\sixbb\scriptscriptfont\bbfam=\fivebb
  \def\bb{\fam\bbfam\ninebb}%
  \textfont\msafam=\ninemsa\scriptfont\msafam=\sixmsa\scriptscriptfont\msafam=\fivemsa
  \textfont\scalnfam=\ninescaln
  \def\scaln{\fam\scalnfam\ninescaln}
  \textfont\ttfam=\ninett
  \def\tt{\fam\ttfam\ninett}%
  \textfont\bffam=\ninebf\scriptfont\bffam=\sixbf\scriptscriptfont\bffam=\fivebf
  \def\bf{\fam\bffam\ninebf}%
  \abovedisplayskip=3pt plus2pt minus 2pt
  \belowdisplayskip=3pt plus1pt minus 2pt
  \abovedisplayshortskip= 2pt plus 2pt minus 1pt
  \belowdisplayshortskip= 2pt plus 2pt minus 1pt
  \smallskipamount=2pt plus 1pt minus 2pt
  \medskipamount=3pt plus 2pt minus 2pt
  \bigskipamount=7pt plus 3pt minus 3pt
  \setbox\strutbox=\hbox{\vrule height 5pt depth 2pt width 0pt}%
  \normalbaselineskip=10.5pt plus.3pt minus.3pt\normalbaselines\rm}

  \def\sevenpoint{%
  \textfont0=\sevenrm\scriptfont0=\sixrm\scriptscriptfont0=\fiverm
  \def\rm{\fam0\sevenrm}%
  \textfont1=\seveni\scriptfont1=\sixi
  \scriptscriptfont1=\fivei\def\oldstyle{\fam1\seveni}%
  \textfont2=\sevensy\scriptfont2=\sixsy\scriptscriptfont2=\fivesy
  \textfont3=\sevenex\scriptfont3=\fiveex
  \textfont\itfam=\sevenit
  \def\it{\fam\itfam\sevenit}%
  \textfont\slfam=\sevensl
  \def\sl{\fam\slfam\sevensl}%
  \textfont\bbfam=\sevenbb \scriptfont\bbfam=\sixbb\scriptscriptfont\bbfam=\fivebb
  \def\bb{\fam\bbfam\sevenbb}%
  \textfont\msafam=\sevenmsa\scriptfont\msafam=\sixmsa
  \textfont\scalnfam=\sevenscaln
  \def\scaln{\fam\scalnfam\sevenscaln}
  \textfont\bffam=\sevenbf\scriptfont\bffam=\sixbf\scriptscriptfont\bffam=\fivebf
  \def\bf{\fam\bffam\sevenbf}%
  \textfont\ttfam=\seventt
  \abovedisplayskip=2pt plus2pt minus 2pt
  \belowdisplayskip=2pt plus1pt minus 2pt
  \abovedisplayshortskip= 1pt plus 2pt minus 1pt
  \belowdisplayshortskip= 1pt plus 2pt minus 1pt
  \smallskipamount=2pt plus 1pt minus 2pt
  \medskipamount=3pt plus 2pt minus 2pt
  \bigskipamount=7pt plus 3pt minus 3pt
  \setbox\strutbox=\hbox{\vrule height 5pt depth 2pt width 0pt}%
  \normalbaselineskip=9pt\normalbaselines\rm}

 \def\onzepts{%
 \textfont0=\elevenrm\scriptfont0=\ninerm
 \textfont1=\elevenib\scriptfont1=\ninei}

\def\douzepts{%
  \textfont0=\twelverm\scriptfont0=\tenrm\def\rm{\fam0\twelverm}%
  \textfont1=\twelveib\scriptfont1=\teni%
  \textfont2=\twelvesy\scriptfont2=\tensy\scriptscriptfont2=\eightsy
  \textfont3=\twelveex
  \textfont\itfam=\twelveti
  \def\it{\fam\itfam\twelveti}%
  \textfont\bffam=\twelvebf\scriptfont\bffam=\tenbf\scriptscriptfont\bffam=\eightbf
  \def\bf{\fam\bffam\twelvebf}%
  \textfont\bbfam=\twelvebb \scriptfont\bbfam=\tenbb
  \def\bb{\fam\bbfam\twelvebb}%
  \textfont\msafam=\twelvemsa\scriptfont\msafam=\tenmsa
  \textfont\scalnfam=\twelvescaln
  \normalbaselineskip=15pt\normalbaselines\rm}

\def\quatorzepts{%
  \textfont0=\fourteenrm\scriptfont0=\twelverm\def\rm{\fam0\fourteenrm}%
  \textfont1=\fourteenib\scriptfont1=\twelveib%
  \textfont2=\fourteensy\scriptfont2=\twelvesy\scriptscriptfont2=\tensy
  \textfont3=\fourteenex
  \textfont\itfam=\fourteenti
  \def\it{\fam\itfam\fourteenti}%
  \textfont\bffam=\fourteenbf\scriptfont\bffam=\twelvebf\scriptscriptfont\bffam=\tenbf
  \def\bf{\fam\bffam\fourteenbf}%
  \textfont\bbfam=\fourteenbb \scriptfont\bbfam=\twelvebb
  \def\bb{\fam\bbfam\fourteenbb}%
  \textfont\msafam=\fourteenmsa\scriptfont\msafam=\twelvemsa
  \textfont\scalnfam=\twelvescaln
  \normalbaselineskip=18pt\normalbaselines\rm}


\def\AA{{\it Acta Arith.}}

\def\picture #1 by #2 (#3){\leavevmode\vbox to #2{
     \hrule width #1 height 0pt depth 0pt
      \vfill
       \special{picture #3}}}

\def\illustration #1 by #2 (#3) scaled #4{\dimen1=#2
  \divide\dimen1 by 1000
  \multiply\dimen1 by #4
  \vtop to \dimen1{\dimen1=#1
  \divide\dimen1 by 1000
  \multiply\dimen1 by #4
  \hsize=\dimen1\vss
  \noindent\special{illustration #3 scaled #4}}}

%% file: option_keys.tex
\ifx\optionkeymacros\undefined\else \fi

\catcode`\Œ=\active\defŒ{{\aa}}       
\catcode`\º=\active\defº{\int}        
\catcode`\=\active\def{\c c}        
\catcode`\¶=\active\def¶{\partial}    
\catcode`\Ä=\active\defÄ{\oint}       
\catcode`\Æ=\active\defÆ{\triangle}   
\catcode`\Â=\active\defÂ{\neg}        
\catcode`\µ=\active\defµ{\mu}         
\catcode`\¿=\active\def¿{{\o}}        
\catcode`\¹=\active\def¹{\pi}         
\catcode`\Ï=\active\defÏ{{\oe}}       
\catcode`\§=\active\def§{{\ss}}       
\catcode`\ =\active\def {\dagger}     
\catcode`\Ã=\active\defÃ{\sqrt}       
\catcode`\·=\active\def·{\Sigma}      
\catcode`\Å=\active\defÅ{\approx}     
\catcode`\½=\active\def½{\Omega}      
\catcode`\£=\active\def£{{\it\$}}     
\catcode`\°=\active\def°{\infty}      
\catcode`\¤=\active\def¤{{\S}}        
\catcode`\¦=\active\def¦{{\P}}        
\catcode`\¥=\active\def¥{\bullet}     
\catcode`\»=\active\def»{\leavevmode\raise.585ex\hbox{\b a}}      
\catcode`\¼=\active\def¼{\leavevmode\raise.6ex\hbox{\b o}}        
\catcode`\­=\active\def­{\not=}       
\catcode`\²=\active\def²{\leq}        
\catcode`\³=\active\def³{\geq}        
\catcode`\Ö=\active\defÖ{\div}        
\catcode`\É=\active\defÉ{{\dots}}     
\catcode`\¾=\active\def¾{{\ae}}       
\catcode`\Ç=\active\defÇ{\og}         
\catcode`\Ò=\active\defÒ{``}          
\catcode`\Á=\active\defÁ{!`}          
\catcode`\¢=\active\def¢{\rlap/c}     
\catcode`\Ô=\active\defÔ{`}           
\catcode`\Õ=\active\defÕ{'}           


\catcode`\=\active\def{{\AA}}       
\catcode`\'=\active\def'{\c C}        
\catcode`\¯=\active\def¯{{\O}}        
\catcode`\¸=\active\def¸{\Pi}         
\catcode`\Î=\active\defÎ{{\OE}}       
\catcode`\®=\active\def®{{\AE}}       
\catcode`\×=\active\def×{\diamond}    
\catcode`\¡=\active\def¡{\accent'27}  
\catcode`\Ó=\active\defÓ{''}          
\catcode`\±=\active\def±{\pm}         
\catcode`\È=\active\defÈ{\fg}         
\catcode`\À=\active\defÀ{?`}          
\catcode`\Ð=\active\defÐ{--}          
\catcode`\Ñ=\active\defÑ{---}         


\catcode`\Š=\active\defŠ{\"a}        
\catcode`\'=\active\def'{\"e}        
\catcode`\•=\active\def•{\"{\i}}     
\catcode`\š=\active\defš{\"o}        
\catcode`\Ÿ=\active\defŸ{\"u}        
\catcode`\Ø=\active\defØ{\"y}        
\catcode`\å=\active\defå{\^A}        
\catcode`\€=\active\def€{\"A}        
\catcode`\…=\active\def…{\"O}        
\catcode`\†=\active\def†{\"U}        
\catcode`\‡=\active\def‡{\'a}        
\catcode`\Ž=\active\defŽ{\'e}        
\catcode`\'=\active\def'{\'{\i}}     
\catcode`\—=\active\def—{\'o}        
\catcode`\œ=\active\defœ{\'u}        
\catcode`\ƒ=\active\defƒ{\'E}        
\catcode`\æ=\active\defæ{\^E}        
\catcode`\é=\active\defé{\`E}        %
\catcode`\ˆ=\active\defˆ{\`a}        
\catcode`\=\active\def{\`e}        
\catcode`\"=\active\def"{\`{\i}}     
\catcode`\˜=\active\def˜{\`o}        
\catcode`\=\active\def{\`u}        
\catcode`\Ë=\active\defË{\`A}        
\catcode`\‹=\active\def‹{\~a}        
\catcode`\–=\active\def–{\~n}        
\catcode`\›=\active\def›{\~o}        
\catcode`\Ì=\active\defÌ{\~A}        
\catcode`\"=\active\def"{\~N}        
\catcode`\Í=\active\defÍ{\~O}        
\catcode`\‰=\active\def‰{\^a}        
\catcode`\=\active\def{\^e}        
\catcode`\"=\active\def"{\^{\i}}     
\catcode`\™=\active\def™{\^o}        
\catcode`\ž=\active\defž{\^u}        

\let\optionkeymacros\null